\documentclass[reqno]{amsart}

\usepackage{amssymb}
\usepackage{amsmath}
\usepackage{setspace}
\usepackage{MnSymbol}

\newtheorem{theorem}{Theorem}[section]
\newtheorem{definition}[theorem]{Definition}
\newtheorem{proposition}[theorem]{Proposition}
\newtheorem{corollary}[theorem]{Corollary}

\newtheorem{remark}[theorem]{Remark}
\newtheorem{example}[theorem]{Example}

\begin{document}

\title{Exponential dichotomies of evolution operators in Banach spaces}

\author{Nicolae Lupa}
\address{N. Lupa: Faculty of Economics and Business Administration,
              West University of Timi\c soara,
              Blvd. Pestalozzi 16, 300115 Timi\c soara, Romania}
\email{nlupa@math.uvt.ro}

\author{Mihail Megan}
\address{M. Megan: Academy of Romanian Scientists,
              Independen\c tei 54, 050094 Bucharest, Romania}
\email{megan@math.uvt.ro}

\date{\today}

\begin{abstract}
This paper considers three dichotomy concepts (exponential dichotomy, uniform exponential dichotomy and strong exponential dichotomy) in the general context of non-invertible evolution operators in Banach spaces. Connections between these concepts are illustrated. Using the notion of Green function, we
give necessary conditions and sufficient ones for strong exponential dichotomy.  Some illustrative examples are presented to prove that the converse of some implication type theorems are not valid.
\end{abstract}

\keywords{Evolution operator, Exponential dichotomy, Lyapunov function}
\subjclass[2000]{47D06, 34D09}

\maketitle

\section{Introduction}

The notion of exponential dichotomy introduced by Perron \cite{Pe.1930} plays a central role in the stability theory of differential equations, discrete dynamical systems \cite{Po.2010}, delay evolution equations \cite{Ca.Re.Sh.2007}, dynamical equations on time scales \cite{Zh.2010}, impulsive equations \cite{Bai.1989}, stochastic processes \cite{St.2010} and many other domains.

The exponential dichotomy property for linear differential equations has gained prominence since the appearance of two fundamental monographs due to Dalecki\v{\i} and Kre\v{\i}n \cite{Da.Kr.1974} and Massera and Sch\"{a}ffer \cite{Ma.Sc.1966}. These were followed by the book of Coppel \cite{Co.1978}, who synthesized and improved the results that existed in the literature up to 1978.
Numerous important papers on this subject appeared afterward and  we  mention in particular \cite{Hu.2006,Mi.Ra.Sc.1998,Po.2006,Pr.Me.1985,Pr.Po.Pr,Sa.2006-1,Sa.2010-1}.
We also refer to the book of Chicone and Latushkin
\cite{Ch.La.1999} for  important results in infinite-dimensional spaces.
In \cite{Sa.Se.1991}, Sacker and Sell use a concept of exponential dichotomy for skew-product semiflows with the restriction that the unstable subspace is finite dimensional. Chow and  Leiva introduce in \cite{Ch.Le.1996-1} a general concept of exponential dichotomy for linear skew-product semiflows weaker that the one used by Sacker and Sell.

One of the most important results in the stability theory of evolution operators is due to Datko \cite{Da.1972} which has given an integral characterization of uniform exponential stability.  This characterization is used to obtain a necessary and sufficient condition  in terms of Lyapunov functions. Preda and Megan extend Datko's result to uniform exponential dichotomy \cite{Pr.Me.1985}.
The importance of Lyapunov functions is well-known in the study of exponential behavior of solutions of differential equations both in finite and infinite-dimensional settings (we refer to the books  \cite{Ch.La.1999,Co.1978,Da.Kr.1974,Ma.Ma.2009,Mi.Sa.Ku.2003}).
Other important Lyapunov type characterizations were obtained in \cite{Ca.Re.Sh.2007,Ha.2001,Me.Bu.1992}.

In the non-autonomous setting, the concept of uniform exponential dichotomy is too restrictive and it is important to consider more general behavior,  for example the nonuniform case, where a consistent contribution is
due to Barreira and Valls \cite{Ba.Va.2008-2,Ba.Va.2009-2,Ba.Va.2009-3,Ba.Va.2011}. Their study is motivated by ergodic  theory and nonuniform hyperbolic theory (we refer the reader to  \cite{Ba.Pe.2002} for details and further information).
Preda and Megan introduced in \cite{Pr.Me.1983} a concept of nonuniform exponential dichotomy which  does not require anything about the norms of the dichotomy projections. In  \cite{Me.Sa.Sa.2002}, the authors involve a  more general type of exponential dichotomy, but with the same assumption on the associated projections. In the more recent papers \cite{Ba.Va.2009-2,Ba.Va.2009-3}, Barreira and Valls  explicitly construct  Lyapunov  functions for given dynamics in the  nonuniform case.

In this paper we consider three dichotomy concepts (exponential dichotomy and uniform exponential dichotomy in Definition \ref{n.e.d} and strong exponential dichotomy in Definition \ref{d4}) in the general context of evolution operators in Banach spaces and we study the connections between these concepts. Some illustrative examples are given to better understand the mathematical setting.

The main goal of this paper is to extend Datko's results to strong  exponential dichotomy.  Thus, we give necessary conditions and sufficient ones for the existence of strong exponential dichotomy in terms of integral inequalities (Proposition \ref{p1} and Theorem\ref{p2}) and Lyapunov functions (Theorem \ref{t.Ly.nec} and Theorem \ref{t.Ly.suf}). As consequences of the above mentioned results, we present necessary and sufficient conditions for uniform exponential dichotomy.
Moreover, we note that we do not need  to assume the invertibility of the evolution operators on the whole  space $X$ (unlike the case of evolution operators generated
by differential equations).


\section{Exponential dichotomies of evolution operators}

Let $X$ be a real or complex Banach space, $\mathcal{B}(X)$ be the Banach algebra of all
bounded linear operators on $X$ and $J$ be a real interval which is $\mathbb{R}_{+}$ or $\mathbb{R}$.  The norms on $X$ and on $\mathcal{B}(X)$  will be denoted by $\parallel \cdot \parallel $. Also, we consider
$$\Delta_{J}^{+}=\left\{ (t,s)\in J\times J: t\geq s \right\}.$$
We first recall the definition of an evolution operator.
\begin{definition}\label{d1} \rm
An operator valued function  $U:\Delta_{J}^{+}\rightarrow\mathcal{B}(X)$ is said to be  an  \emph{evolution operator} (on $J$) if the following conditions are satisfied:

\begin{enumerate}
\item $U(t,t)=Id$ (the identity operator on $X$) for every $t\in J$;

\item $U(t,s)U(s,t_0 )=U(t,t_0)$, for all $ t\geq s\geq t_0$ in $J$;

\item $\Delta_J^+\ni(t,s)\longmapsto U(t,s)x\in X$ is continuous for every $x\in X$.
\end{enumerate}
\end{definition}
The concept of evolution operator arises naturally from the theory of well-possed evolution equations.
Roughly speaking, when the Cauchy problem
\begin{equation}\label{eq.Cauchy}
\quad\left\{\begin{array}{ll}
\dot{u}(t)=A(t)u(t),\;t\geq{s}\in J\\
u(s)=x
\end{array}\right.
\end{equation}
is well-posed with regularity subspaces $\left(Y_t\right)_{t\in J}$, then the operator
$$U(t,s)x:= u(t;s,x) \text{ for } t\geq s \text{ and } x\in Y_s,$$
where $u(\cdot;s,x)$ is the unique solution of Eq. (\ref{eq.Cauchy}), can be extended by continuity to an evolution operator.
For more details on well-posed non-autonomous Cauchy problems  we refer the reader to Nagel and Nickel \cite{Na.Ni.2002} and the references therein.

\begin{definition}\label{d2}\rm
A strongly continuous  function $P:J\rightarrow \mathcal{B}(X)$ is said to be a \emph{projection valued function} if $P^{2}(t)=P(t),\text{ for every } t\in J.$  

If $P:J\rightarrow \mathcal{B}(X)$ is a projection valued function, we denote by $Q(t)=Id-P(t)$, the complementary projection of $P(t)$ for each $t\in J$.
One can easily see that
\begin{equation*}\label{eq.projection}
P(t)Q(t)=Q(t)P(t)=0 \text{ and } \text{Ker}\,P(t)=\text{Range}\,Q(t),\text{ for every } t\in J.
\end{equation*}
\end{definition}

\begin{definition}\label{d3}\rm
Given an evolution operator $U:\Delta_{J}^{+}\rightarrow\mathcal{B}(X)$, we say that a projection valued function  $P:J\rightarrow \mathcal{B}(X)$ is \emph{compatible with ${U}$} if the following conditions hold:
\begin{enumerate}
\item $P(t)U(t,s)=U(t,s)P(s), \text{ for all }(t,s)\in\Delta_{J}^{+}$;

\item The restriction $U(t,s)_{| \text{Ker}\, P(s)}:\text{Ker}\, P(s)\rightarrow \text{Ker}\, P(t)$ is an isomorphism for each $(t,s)\in\Delta_{J}^{+}$; we denote its inverse by $V(s,t)$.
\end{enumerate}
\end{definition}

\begin{proposition}
If $P:J\rightarrow \mathcal{B}(X)$ is a projection valued function compatible with the evolution operator $U:\Delta_{J}^{+}\rightarrow\mathcal{B}(X)$, then for any $(t,s)\in\Delta_{J}^{+}$, the operator $V(s,t)$ is an isomorphism from Range $Q(t)$ to Range $Q(s)$  and the following properties hold:
\begin{enumerate}
\item[$(v_1)$] $U(t,s)V(s,t)Q(t)=Q(t)$; \label{v1}
\item[$(v_2)$] $V(s,t)U(t,s)Q(s)=Q(s)$;  \label{v2}
\item[$(v_3)$] $V(s,t)Q(t)=Q(s)V(s,t)Q(t)$.  \label{v3}
\end{enumerate}
Moreover, we have that
\begin{enumerate}
\item[$(v_4)$] $V(t_0,s)V(s,t)Q(t)=V(t_0,t)Q(t)$, for all $t\geq s\geq t_0$ in $J$. \label{v4}
\end{enumerate}
\end{proposition}
\begin{proof}
Indeed, relations $(v_1)$--$(v_3)$ are simple consequences of the definition of the operator $V(s,t)$. In order to prove $(v_4)$ we use $(v_1)$--$(v_3)$ and Definition \ref{d1} and  obtain
\begin{align*}
V(t_0,s)V(s,t)Q(t)&=V(t_0,s)V(s,t)U(t,t_0)V(t_0,t)Q(t)\\
&=V(t_0,s)V(s,t)U(t,s)U(s,t_0)Q(t_0)V(t_0,t)Q(t)\\
&=V(t_0,s)V(s,t)U(t,s)Q(s)U(s,t_0)V(t_0,t)Q(t)\\
&=V(t_0,s)Q(s)U(s,t_0)V(t_0,t)Q(t)\\
&=V(t_0,s)U(s,t_0)Q(t_0)V(t_0,t)Q(t)\\
&=Q(t_0)V(t_0,t)Q(t)\\
&=V(t_0,t)Q(t).
\end{align*}
\end{proof}

\begin{example}\label{ex.baza}\rm
Let $X$ be a real Banach algebra with the unit ${\bf{1}}$ and $\parallel {\bf{1}} \parallel=1$. We denote by $l^{\,\infty}(\mathbb{N},X)$ the space of all bounded sequences $x=(x_0,x_1,x_2,\ldots, x_n,\ldots)$ with $x_n\in X$, $n\in\mathbb{N}$. $l^{\,\infty}(\mathbb{N},X)$ is a real Banach space endowed with the natural norm
$$\| x \|_{\infty}=\sup\limits_{n\in\mathbb{N}} \parallel x_n \parallel.$$\\
Let $u,v:J\rightarrow (0,\infty)$ be two continuous functions and consider the operator $$U(t,s):l^{\,\infty}(\mathbb{N},X)\rightarrow l^{\,\infty}(\mathbb{N},X),$$ defined by
\begin{equation}\label{eq.op.ev.1}
U(t,s)(x_0,x_1,x_2,x_3,\ldots)=\left(\frac{u(s)}{u(t)}x_0,\frac{v(t)}{v(s)}x_1,\frac{u(s)}{u(t)}x_2,\frac{v(t)}{v(s)}x_3,\ldots\right),\text{ for  $(t,s)\in\Delta_{J}^{+}.$}
\end{equation}
Let now $P:J\rightarrow \mathcal{B}\left(l^{\,\infty}(\mathbb{N},X)\right)$ be the operator valued function defined by
\begin{equation}\label{eq.projector}
P(t)(x_0,x_1,x_2,x_3,\ldots)=(x_0,0,x_2,0,\ldots), \text{ for } t\in J.
\end{equation}
It is  easy to see that $U$ is an evolution operator on  $l^{\,\infty}(\mathbb{N},X)$ and $P(\cdot)$ is a projection valued function with the complementary projection
$$Q(t)(x_0,x_1,x_2,x_3,\ldots)=(0,x_1,0,x_3,\ldots), \text{ for } t\in J.$$
We observe that
\begin{equation}\label{eq.op.ev}
U(t,s)=\frac{u(s)}{u(t)} P(s)+\frac{v(t)}{v(s)} Q(t), \text{ for } (t,s)\in \Delta_{J}^{+}.
\end{equation}
Since Range$P(t)$=Range $P(s)$ for all $t,s\in J$, it follows that
\begin{equation}\label{eq.comut}
P(t)P(s)=P(s) \text{ and }P(s)P(t)=P(t), \text{ for }(t,s)\in \Delta_{J}^{+}.
\end{equation}
Combining (\ref{eq.op.ev}) and (\ref{eq.comut}) we obtain  that  $P(\cdot)$ is  compatible with $U$. Furthermore, we have
\begin{equation}\label{eq.per}
U(t,s)P(s)= \frac{u(s)}{u(t)} P(s) \text{ and } V(s,t)Q(t)= \frac{v(s)}{v(t)} Q(s), \text{ for } (t,s)\in \Delta_{J}^{+}.
\end{equation}
\end{example}

In their notable work \cite{Ba.Va.2008-2,Ba.Va.2009-2,Ba.Va.2009-3,Ba.Va.2011}, Barreira and Valls consider the following notion of (nonuniform) exponential dichotomy:

\begin{definition}\label{n.e.d}\rm
An evolution operator $U:\Delta_{J}^{+}\rightarrow\mathcal{B}(X)$ is said to have an \emph{exponential dichotomy} (in $J$) if there is
a projection valued function $P:J\rightarrow \mathcal{B}(X)$ compatible with $U$ and there exist constants $N\geq 1$, $\alpha\geq 0$ and $\beta>0$ such that
\begin{equation*}
\parallel U(t,s)P(s)\parallel\leq Ne^{\alpha |s|}e^{-\beta(t-s)} \text{ and }
\parallel V(s,t)Q(t)\parallel \leq Ne^{\alpha |t|}e^{-\beta(t-s)},
\end{equation*}
for all $(t,s)\in\Delta_{J}^{+}$.

The constant $\alpha$ measures the non-uniformity of the dichotomy. In particular, when $\alpha=0$ the evolution operator  $U$ is said to have a \emph{uniform exponential dichotomy}. The projection valued function $P(\cdot)$ will be called the \emph{dichotomy projection}.
\end{definition}

\begin{remark}\label{rem.eq}\rm
An evolution operator $U:\Delta_{J}^{+}\rightarrow\mathcal{B}(X)$ has an \emph{exponential dichotomy} if and only if there exist
a projection valued function $P:J\rightarrow \mathcal{B}(X)$ compatible with $U$ and constants $N\geq 1$, $\alpha\geq 0$ and $\beta>0$ such that
\begin{enumerate}

\item[(ed$_1$)] $\parallel U(t,s)P(s)x\parallel\leq Ne^{\alpha |s|}e^{-\beta(t-s)}\parallel P(s)x\parallel$, for $(t,s,x)\in\Delta_{J}^{+}\times X$;

\item[(ed$_2$)] \;\;$\;Ne^{\alpha |t|}\parallel U(t,s)Q(s)x\parallel \geq e^{\beta(t-s)}\parallel Q(s)x\parallel$, for $(t,s,x)\in\Delta_{J}^{+}\times X$;

\item[(ed$_3$)] \;\;$\parallel P(t)\parallel\leq N e^{\alpha |t|}$, for $t\in J$.

\end{enumerate}
\end{remark}
\begin{proof}
Indeed, if we assume that $U$  has an  exponential dichotomy with the dichotomy projection $P(\cdot)$ and constants $N\geq 1$, $\alpha\geq 0$ and $\beta>0$, we have
\begin{align*}
\parallel U(t,s)P(s)x\parallel=\parallel U(t,s)P(s)P(s)x\parallel\leq Ne^{\alpha |s|}e^{-\beta(t-s)}\parallel P(s)x\parallel
\end{align*}
and, respectively
\begin{align*}
\parallel Q(s)x\parallel = \parallel V(s,t)Q(t)U(t,s)Q(s)x\parallel\leq Ne^{\alpha |t|}e^{-\beta(t-s)}\parallel U(t,s)Q(s)x\parallel,
\end{align*}
for all $(t,s,x)\in\Delta_{J}^{+}\times X$. Setting $t=s$ in Definition  \ref{n.e.d} we obtain (ed$_3$).

Conversely, for $x\in X$ with $\parallel x\parallel=1$ and $(t,s)\in\Delta_{J}^{+}$, successively we have
\begin{align*}
\parallel U(t,s)P(s)x\parallel\leq N e^{\alpha |s|} e^{-\beta(t-s)}\parallel P(s)x\parallel\leq N^2 e^{2\alpha |s|} e^{-\beta(t-s)}
\end{align*}
and
\begin{align*}
\parallel V(s,t)Q(t)x\parallel &=\parallel Q(s) V(s,t)Q(t)x\parallel\\
&\leq Ne^{\alpha |t|}e^{-\beta(t-s)}\parallel U(t,s)Q(s)V(s,t)Q(t)x\parallel\\
&=Ne^{\alpha |t|}e^{-\beta(t-s)}\parallel Q(t)x\parallel\leq Ne^{\alpha |t|}e^{-\beta(t-s)}(1+\parallel P(t)\parallel)\\
&\leq 2N^2 e^{2\alpha |t|} e^{-\beta(t-s)}.
\end{align*}
Hence, $U$ has an exponential dichotomy.
\end{proof}

We note that relation (ed$_2$) implies that the operator $U(t,s)$ is injective from the range of $Q(s)$ to the range of $Q(t)$.
In the uniform case, if the evolution operator  $U$  is exponentially bounded (i.e. there exist constants $M\geq 1$ and $\omega>0$ such that $|| U(t,s) ||\leq M e^{\omega(t-s)}$ for all $(t,s)\in\Delta_{J}^{+}$), relation (ed$_3$) can be omitted (see Lemma 4.2 in \cite{Mi.Ra.Sc.1998}).

Geometrically, if an evolution operator on the real line ($J=\mathbb{R}$) is generated by a well-posed evolution equation, an exponential  dichotomy (uniform or not)  means that the space $\mathbb{R}\times X$ splits  into two invariant vector bundles consisting of solutions with a specific  asymptotic behavior (see \cite{Po.2010}), namely
\begin{itemize}
\item the stable bundle ($X_s(s)$) consisting of exponentially forward solutions on the set $\Delta^{+}_{s}=\{t\in \mathbb{R}:t\geq s\}$ and given by the range of $P(s)$, where $P(\cdot)$ is the dichotomy projection;
\item the unstable bundle ($X_u(s)$) consisting of solutions which exist in backward time on $\Delta^{-}_{s}=\{t\in \mathbb{R}:t\leq s\}$ and are exponentially decaying, given by the range of $Q(s)$.
\end{itemize}

In the general case of an abstract evolution operator  $U$ which has a  certain type of exponential dichotomy with the dichotomy projection $P(\cdot)$, the previous statements can be translated as:
\begin{equation}\label{rel.stable}
U(t,s)P(s)x\rightarrow 0 \text{ (exponentially) as } t\rightarrow+\infty
\end{equation}
and
\begin{equation}\label{rel.unstable}
V(t,s)Q(s)x\rightarrow 0 \text{ (exponentially) as } t\rightarrow-\infty
\end{equation}
for each $(s,x)\in\mathbb{R}\times X$.
For example, if an evolution operator $U:\Delta_{\mathbb{R}}^{+}\rightarrow\mathcal{B}(X)$ has an exponential dichotomy with the dichotomy projections $P(\cdot)$ (as in Definition \ref{n.e.d}), then  relations (\ref{rel.stable}) and (\ref{rel.unstable}) hold.

Otherwise, if  $J=\mathbb{R}_{+}$, relation (\ref{rel.unstable}) must be replaced with a proper one such that we can characterize the unstable bundle in terms of asymptotic behavior of  trajectories.
From Remark \ref{rem.eq}, we observe that
if an evolution operator $U:\Delta_{\mathbb{R}_{+}}^{+}\rightarrow\mathcal{B}(X)$  has a uniform exponential dichotomy with the dichotomy projection $P:\mathbb{R}_{+}\rightarrow \mathcal{B}(X)$, then for any   $s\geq 0$ and $x\in X$, relation (\ref{rel.stable}) holds and we have that
\begin{equation}\label{rel.2}
\parallel U(t,s)Q(s)x\parallel\rightarrow \infty \text{ (exponentially) as } t\rightarrow\infty,  \text{ unless } Q(s)x=0.
\end{equation}

In the following we only consider evolution operators on the half-line ($J=\mathbb{R}_{+}$) and we denote by $\Delta_{+}$ the set $\Delta_{\mathbb{R}_{+}}^{+}$.

\begin{definition}\label{d4}\rm
We say that an evolution operator $U:\Delta_{+}\rightarrow\mathcal{B}(X)$ has a \emph{strong exponential
dichotomy} if there is a projection valued function $P:\mathbb{R}_{+}\rightarrow \mathcal{B}(X)$  compatible with $U$ and there exist constants
$N\geq 1$, $\alpha\geq 0$ and $\beta>0$ such that
\begin{equation*}
\parallel U(t,s)P(s)\parallel\leq Ne^{\alpha s}e^{-\beta(t-s)}  \text{ and } \parallel V(s,t)Q(t)\parallel\leq Ne^{\alpha s}e^{-\beta(t-s)},
\end{equation*}
for all $(t,s)\in\Delta_{+}$.
\end{definition}

\begin{remark}\label{rem.pr}\rm
If the evolution operator  $U$ has a strong exponential
dichotomy then there exist a projection valued function $P(\cdot)$ compatible with $U$ and constants $N\geq 1$, $\alpha\geq 0$ and $\beta>0$ such that for all $(t,s)\in\Delta_{+}$ and $x\in X$ we have
$$\parallel U(t,s)P(s)x\parallel\leq Ne^{\alpha s}e^{-\beta(t-s)}\parallel P(s) x\parallel$$ { and }
$$\;\;N e^{\alpha s} \parallel U(t,s)Q(s)x\parallel \geq e^{\beta(t-s)}\parallel Q(s)x\parallel.$$
The inequalities above prove that relations \eqref{rel.stable} and \eqref{rel.2} are fulfilled.
Moreover, the dichotomy projection is exponentially bounded, i.e.  $\parallel P(t)\parallel\leq N e^{\alpha t}$, for every $t\geq 0$.
\end{remark}

\begin{proposition}\label{rem.equivalence}
 An evolution operator  $U$ has a strong exponential dichotomy if and only if there is a projection valued function $P(\cdot)$ compatible with $U$ and there exist constants
$N_{i}\geq 1$, $\alpha_{i}\geq 0$ and $\beta_{i}>0$, $i=1,2$ with $\alpha_{2}<\beta_{2}$ such that
\begin{equation*}
\parallel U(t,s)P(s)\parallel\leq N_{1}e^{\alpha_{1} s} e^{-\beta_{1}(t-s)} \text{ and }
\parallel V(s,t)Q(t)\parallel\leq N_{2}e^{\alpha_{2} t} e^{-\beta_{2}(t-s)}, \text{ for  $(t,s)\in\Delta_{+}$.}
\end{equation*}
\end{proposition}
\begin{proof} It is a simple exercise. Indeed, if $U$ has a strong exponential dichotomy with the dichotomy projection $P(\cdot)$ and constants
$N\geq 1$, $\alpha\geq 0$ and $\beta>0$, we have
\begin{align*}
\parallel V(s,t)Q(t)\parallel\leq Ne^{\alpha s}e^{-\beta(t-s)}=Ne^{\alpha t}e^{-(\alpha+\beta)(t-s)}, \text{ for } (t,s)\in\Delta_{+}.
\end{align*}
Conversely, in a similar manner, it results that $U$ has a strong exponential dichotomy with constants $N=\max\{N_1, N_2\}$, $\alpha=\max\{\alpha_1,\alpha_2\}$ and $\beta=\min\{\beta_1,\beta_2-\alpha_2\}$.
\end{proof}

We remark that due to the invertibility assumption in Definition \ref{d3}, the second inequality in the proposition above is just a version of the first one when time goes backwards (like in Definition \ref{n.e.d}). Additionally,  we impose that the nonuniform part on the unstable bundle is smaller than the uniform one. However, this is not a restrictive assumption since the notion of (nonuniform) exponential dichotomy is apparently motivated by ergodic theory (see \cite{Ba.Pe.2002} and the references therein) which states that the nonuniform part in the dichotomies of the ``most'' equations is arbitrarily small (is as small as desired in comparison to the Lyapunov exponents). The concept of strong exponential dichotomy is also motivated by Theorem 2 in \cite{Ba.Va.2008-2} which states that if a non-autonomous linear equation has a strong exponential dichotomy then, under certain conditions, the perturbed equation has an exponential dichotomy. We notice that the assumption $\alpha<\beta$ in the above mentioned theorem is only used for the unstable part.

It is clear that the existence of (strong) exponential dichotomy does not imply the same type of exponential growth on the whole space $X$. Therefore, it is natural to consider the following concept of exponential growth for a given evolution operator $U$ and a projection valued function $P(\cdot)$ which is compatible with $U$:
\begin{definition}\label{e.g}\rm
Let $U:\Delta_{+}\rightarrow\mathcal{B}(X)$ be an evolution operator and $P:\mathbb{R}_{+}\rightarrow \mathcal{B}(X)$ be a projection valued function compatible with $U$. We say that $U$ has  \emph{$P$-exponential growth} if there exist $M\geq 1$, $\varepsilon\geq 0$ and $\omega>0$ such that
\begin{equation*}\label{eq.gr}
\parallel U(t,s)P(s)\parallel\leq Me^{\varepsilon s}\;e^{\omega(t-s)} \text{ and }\parallel V(s,t)Q(t)\parallel\leq Me^{\varepsilon s}\;e^{\omega(t-s)},
\end{equation*}
for all $(t,s)\in\Delta_{+}$. If $\varepsilon=0$, we say that $U$ has  \emph{$P$-uniform exponential growth}.
\end{definition}

Now, we point out the connections between the dichotomy concepts considered above (uniform exponential dichotomy (u.e.d), strong exponential dichotomy (s.e.d) and exponential dichotomy (e.d)).

If an evolution operator $U$  has a uniform exponential dichotomy then it also has a strong exponential
dichotomy. The following example shows that the converse is not true.

\begin{example}\label{ex.2}\rm
Let
$u(t)=v(t)=e^{t(3 +\cos^2{t})}$, for  $t\geq 0$, in Example \ref{ex.baza}.  Then the evolution operator $U$ defined by (\ref{eq.op.ev.1}) has a strong exponential dichotomy with the dichotomy projection $P(\cdot)$ considered in (\ref{eq.projector}) that is not a uniform one.
\end{example}
\begin{proof}
Indeed, since $\|P(s)\|=\|Q(s)\|=1$ for $s\geq 0$, by relation (\ref{eq.per}) we have
$$\parallel U(t,s)P(s)\parallel=e^{-3(t-s)-t\cos^2{t}+s\cos^2{s}}
\leq e^{s}e^{-3(t-s)}$$
and, respectively
$$\parallel V(s,t)Q(t)\parallel
\leq e^{s}e^{-3(t-s)},$$
for all $(t,s)\in \Delta_{+}$. Thus, the evolution operator $U$ has a strong exponential dichotomy.
If we now assume that $U$ has a uniform exponential dichotomy with the dichotomy projection
$P(\cdot)$,
then there exist
$N\geq 1$ and $\beta>0$ such that
 $$\parallel U(t,s)P(s)\parallel\leq Ne^{-\beta(t-s)}, \text{ for all } (t,s)\in \Delta_{+}.$$
In particular, for  $t=n\pi+\pi/2$ and $s=n\pi$ with $n\in \mathbb{N}$, we have
$$e^{-3\pi/2+n\pi}\leq Ne^{-\beta\pi/2},\text{ for all } n\in \mathbb{N}.$$
Letting $n\rightarrow\infty$ in the relation above, we obtain a contradiction.
\end{proof}

If an evolution operator has a strong exponential dichotomy then it also has an exponential
dichotomy.
The example below shows that the existence of  exponential dichotomy with a dichotomy projection $P(\cdot)$ does not imply the strong exponential dichotomy with the same dichotomy projection.

\begin{example}\label{ex.inst}\rm
If we consider $u(t)=e^{t}$ and $v(t)=e^{-t}$, for $t\geq 0$, in Example \ref{ex.baza}, then
the evolution operator $U$ has an exponential dichotomy with the dichotomy projection $P(\cdot)$ that is not a strong exponential dichotomy.
\end{example}
\begin{proof}
Indeed, again by (\ref{eq.per}), we have
$$\parallel U(t,s)P(s)\parallel= e^{-(t-s)} \text{ and } \parallel V(s,t)Q(t)\parallel=e^{t-s}\leq e^{2t}e^{-(t-s)},$$
for all $(t,s)\in \Delta_{+}$. Hence, $U$ has an exponential dichotomy.
Since  $U(t,s)x\rightarrow 0$ as  $t\rightarrow\infty$, for all $s\geq 0$ and $x\in l^{\,\infty}(\mathbb{N},X)$, we have that relation (\ref{rel.2}) is not valid, and thus, $U$ does not have a strong exponential dichotomy with the dichotomy projection $P(\cdot)$ considered in (\ref{eq.projector}).
\end{proof}

\begin{remark}\rm
The connections between the dichotomy concepts considered above can be synthesized in the following diagram:
$$
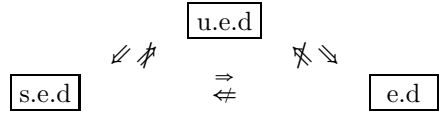
\begin{array}{ccccc}
 &  &\framebox{\makebox{u.e.d}} & &\\
 &\Swarrow \not\Nearrow & & \not\Nwarrow\Searrow   &\\
 \framebox{\makebox{s.e.d}} & & \overset{\Rightarrow}{\nLeftarrow} & & \framebox{\makebox{\,\,e.d\,\,}}\\
\end{array}
$$
We note that the counterexamples between each two concepts are given with the same projection valued function and thus, the diagram is valid
in the case of dichotomy concepts with an a priori defined dichotomy projection.
\end{remark}

In contrast to  the uniform case, where the dichotomy projection $P(\cdot)$ is uniformly bounded, for (strong) exponential dichotomy it is possible that  $|| P(t) ||\rightarrow \infty$ as $t\rightarrow\infty$. Even so, by Remark \ref{rem.pr}, the norms of the dichotomy projections cannot increase faster than an exponential.
\begin{example}\label{ex.proj}\rm
Consider the evolution operator
$$U_{S}(t,s)=S(t)U(t,s)S(s)^{-1}, \text{ for } (t,s)\in \Delta_{+},$$ where $U$ is the evolution operator in Example \ref{ex.2} and $S(t):l^{\,\infty}(\mathbb{N},X)\rightarrow l^{\,\infty}(\mathbb{N},X)$ is defined by
\begin{equation*}
S(t)(x_0,x_1,\ldots)=\left(x_0+\frac{t+1}{\sqrt{1+(t+1)^2}}\,x_1, \frac{1}{\sqrt{1+(t+1)^2}}\,x_1, \ldots\right), \text{ for } t\geq 0.
\end{equation*}
It is easy to see that $S(t)^{-1}$ exits and it is given by
\begin{equation*}
S(t)^{-1}(x_0,x_1,\ldots)=\left(x_0-(t+1)x_1, \sqrt{1+(t+1)^2}x_1, \ldots\right), \text{ for } t\geq 0.
\end{equation*}
For  each $t\geq 0$ and $x=(x_0,x_1,x_2,x_3,\ldots)\in l^{\,\infty}(\mathbb{N},X)$, we take
$$\widetilde P(t)x=(x_0- (t+1)\,x_1,0,x_2- (t+1)\,x_3,0,\ldots).$$
A simple computation shows that $\widetilde P(\cdot)$ is a projection valued function compatible with the evolution operator $U_{S}$ and
\begin{equation*}\label{eq.s1}
\parallel \widetilde{P}(t) x\parallel_{\infty}\leq (2+t)\parallel x\parallel_{\infty}, \text{ for all } t\geq 0 \text{ and } x\in l^{\,\infty}(\mathbb{N},X).
\end{equation*}
Since $\parallel \widetilde{P}(t)(\textbf{1},-\textbf{1},\textbf{1},-\textbf{1},\ldots)\parallel_{\infty}=2+t$, for $t\geq 0$, we deduce that
$$\parallel \widetilde{P}(t) \parallel=2+t\leq e e^{t}, \text{ for all } t\geq 0.$$
Now we observe  that
\begin{align*}
\parallel S(t)x\parallel_{\infty}&=\sup\limits_{n\in\mathbb{N}}\left\{ \parallel x_{2n}+\frac{t+1}{\sqrt{1+(t+1)^2}} \; x_{2n+1}\parallel,\,\frac{1}{\sqrt{1+(t+1)^2}} \parallel x_{2n+1}\parallel\right\}\\
&\leq \sup\limits_{n\in\mathbb{N}}\left\{ \parallel x_{2n}\parallel+ \parallel x_{2n+1}\parallel,\, \parallel x_{2n+1}\parallel
\right\}
\leq 2 \parallel x\parallel_{\infty},
\end{align*}
for all  $t\geq 0$ and $x=(x_0,x_1,x_2,x_3,\ldots)\in l^{\,\infty}(\mathbb{N},X)$ and thus $\parallel S(t)\parallel\leq 2$, for $t\geq 0$.

One can easily verify that $$S(t)^{-1}P(t)=P(t)\text{ and }P(t)\widetilde P(t)=\widetilde P(t), \text{ for }t\geq 0,$$ where $P(\cdot)$ is the projection valued function given by (\ref{eq.projector}).
Since $U$ has a strong exponential dichotomy with the dichotomy projection $P(\cdot)$, it follows that for all $(t,s)\in \Delta_{+}$,  we have
\begin{align*}
\parallel U_{S}(t,s)\widetilde{P}(s)\parallel&=\parallel S(t)U(t,s)S(s)^{-1}P(s)\widetilde{P}(s)\parallel\\
&\leq {2} \parallel U(t,s)P(s)\parallel\parallel\widetilde{P}(s)\parallel\\
&\leq 2ee^{2s}e^{-3(t-s)}.
\end{align*}
On the other hand, since the complementary projection of $\widetilde P(t)$ is given by $$\widetilde{Q}(t)(x_0,x_1,x_2,x_3,\ldots)=((t+1)\,x_1,x_1,(t+1)\,x_3,x_3,\ldots),$$
we deduce that
\begin{align*}
\parallel V_{S}(s,t)\widetilde{Q}(t)x\parallel_{\infty}&= \parallel S(s)V(s,t)S(t)^{-1}((t+1)\,x_1,x_1,(t+1)\,x_3,x_3,\ldots)\parallel_{\infty}\\
&=\parallel S(s)V(s,t)(0,\sqrt{1+(t+1)^2}x_1,0,\sqrt{1+(t+1)^2}x_3,\ldots)\parallel_{\infty}\\
&=\sqrt{1+(t+1)^2} \parallel S(s)V(s,t)Q(t)x\parallel_{\infty}\\
&\leq 2\sqrt{2e}e^{3/2t}e^{-4(t-s)}\parallel x\parallel_{\infty},
\end{align*}
for all  $(t,s)\in \Delta_{+}$ and $x=(x_0,x_1,x_2,x_3,\ldots)\in l^{\,\infty}(\mathbb{N},X)$. The last inequality yields from Example \ref{ex.2}, having in mind that
$$1+\frac{(t+1)^2}{2}\leq e^{t+1}, \text{ for all } t\geq 0.$$
By Proposition \ref{rem.equivalence} it follows that the evolution operator $U_{S}$ has a strong exponential dichotomy with the dichotomy projection $\widetilde P(\cdot)$ and we have
$$\parallel \widetilde P(t)\parallel=2+t\rightarrow\infty \text{ as } t\rightarrow\infty.$$
\end{example}


\section{The main results}

As in \cite{Ba.Va.2011}, we assume the existence of a projection valued function compatible with a
given evolution operator (this assures the existence of a splitting of the space $X$ into two closed subspaces such that the evolution operator leaves invariant the splitting and it is an isomorphism on the unstable part).
Abandoning this assumption causes several technical complications and it is not the purpose of this paper. For instance, this impediment can be eliminated using the notion of admissibility (see \cite{Ch.La.1999,Hu.2006,Pr.Po.Pr,Sa.2006-1}).

For a given evolution operator $U:\Delta_{+}\rightarrow\mathcal{B}(X)$ and projection valued function $P:\mathbb{R}_{+}\rightarrow \mathcal{B}(X)$  compatible with $U$, we denote by
\[
G(t,s)=\begin{cases}
\phantom{-}U(t,s)P(s),&t> s\geq 0\\
-V(t,s)Q(s),&s> t\geq 0
\end{cases}
\]
the Green function of the evolution operator $U$ and the projection valued function $P(\cdot)$ which is compatible with $U$.
Additionally, if $U$ has  $P$-exponential growth, we have that
$$\| G(t,s)\|\leq M e^{\varepsilon s}e^{\omega|t-s|}, \text{ for all } t,s\geq 0,\; t\neq s,$$
where the constants $M\geq 1$, $\varepsilon\geq 0$ and $\omega>0$ are given by Definition \ref{e.g}.


\subsection{A Datko type theorem}
The next result gives a necessary condition for the existence of strong exponential dichotomy:

\begin{proposition}\label{p1}
Let $p>0$ be a real constant. If the evolution operator $U:\Delta_{+}\rightarrow\mathcal{B}(X)$ has a strong exponential
dichotomy with the dichotomy projection $P(\cdot)$ then $U$ has  $P$-exponential growth and there exist  three constants $K\geq 1$, $\delta\geq 0$ and  $\gamma>0$  such that
\begin{equation}\label{eq.Datko}
\int_{0}^{\infty}e^{p\gamma|\tau-t|}\parallel G(\tau,t)x\parallel ^{p}d\tau
\leq Ke^{p\delta t}|| x ||^p, \text{ for all  $t\geq 0$ and $x\in X$.}
\end{equation}
\end{proposition}
\begin{proof}
It is easy to see that $U$ has $P$-exponential growth, taking $M=N$, $\varepsilon=\alpha$ and $\omega>0$, an arbitrary  positive real number.  Furthermore,  relation (\ref{eq.Datko}) holds for $\delta=\alpha$, $\gamma\in(0,\beta)$ and $K=\max\left\{\frac{2N^p}{p(\beta-\gamma)},1\right\}$, where $N\geq 1$, $\alpha\geq 0$ and  $\beta>0$ and  are given by Definition \ref{d4}.
\end{proof}

The example below shows that the converse of the previous result may not be valid:

\begin{example}\label{ex.cond}\rm
Consider the evolution operator $U$  and the projection valued function $P(\cdot)$ in Example \ref{ex.baza}, for  $u(t)=e^{t(1+\cos^2{t})}$ and $v(t)=e^{t\cos^2{t}}$, $t\geq 0$.

$U$ has $P$-exponential growth and for any positive real parameter $p>0$,
we have
\begin{align*}
\int_{0}^{\infty}e^{\frac{1}{2}p|\tau-t|}&\parallel G(\tau,t)x\parallel_{\infty} ^{p}d\tau
\leq\int_{t}^{\infty}e^{-\frac{1}{2}p(\tau-t)}e^{-p\tau\cos^2{\tau}+pt\cos^2{t}}d\tau\|x\|_{\infty}^p+\\
&+\int_{0}^{t}e^{-\frac{1}{2}p(t-\tau)}e^{-p\tau\sin^2{\tau}+p t\sin^2{t}}d\tau \|x\|_{\infty}^p\\
&\leq e^{pt}\left(\int_{t}^{\infty}e^{-\frac{1}{2}p(\tau-t)}d\tau+\int_{0}^{t}e^{-\frac{1}{2}p(t-\tau)}d\tau\right)||x||_\infty^p\\
&\leq (4/p)e^{p t} ||x||_\infty^p, \text{ for  $t\geq 0$ and $x\in l^{\,\infty}(\mathbb{N},X).$}
\end{align*}
Thus,  (\ref{eq.Datko}) holds. If we assume that the evolution operator $U$ has a strong exponential dichotomy with the dichotomy projection $P(\cdot)$, then there
exist  $N\geq 1$, $\alpha\geq 0$ and $\beta>0$ such that
\begin{equation*}
\parallel V(s,t)Q(t)\parallel\leq Ne^{\alpha s}e^{-\beta(t-s)},
\end{equation*}
for all $(t,s)\in\Delta_{+}$. Letting $t=2n\pi+\frac{\pi}{2}$ with $n\in\mathbb{N}$ and $s=0$, we obtain
$$e^{(2n\pi+\pi/2)\beta}\leq N,\text{ for all } n\in\mathbb{N},$$
which is false. Hence, the evolution operator $U$ does not have  a strong exponential dichotomy with  the dichotomy projection
$$P(t)(x_0,x_1,x_2,x_3,\ldots)=(x_0,0,x_2,0,\ldots), \text{ for } t\geq 0.$$
\end{example}

\begin{remark}\rm
The result in Proposition \ref{p1} remains valid even if the evolution operator $U$ has an exponential dichotomy.
\end{remark}

The question is: What additional properties must the constants $\delta$ and $\gamma$  posses such that (\ref{eq.Datko}) implies the existence of strong exponential dichotomy for the evolution operator $U$?

To answer this question we now come  to our first main result, which is a  Datko  type theorem for strong exponential dichotomy. A version of the theorem below is given in \cite{Lu.Me.Po.2010} for weak exponential stability and in \cite{Lu.Me.2010} for nonuniform exponential stability, using different approaches.

\begin{theorem}\label{p2}
Let $U:\Delta_{+}\rightarrow\mathcal{B}(X)$ be an evolution operator and $P:\mathbb{R}_{+}\rightarrow \mathcal{B}(X)$ be a projection valued function compatible with $U$.
If $U$ has $P$-exponential growth and there exist  $p>0$, $K\geq 1$, $\gamma>\varepsilon$ and $\delta\in[0,\gamma)$ such that $(\ref{eq.Datko})$ holds
then $U$ has a strong exponential dichotomy.
\end{theorem}
\begin{proof}
Considering  $(t,s)\in\Delta_{+}$ with $t\geq s+1$ and $x\in X$,  we obtain
\begin{align*}
e^{p\gamma (t-s)}\parallel U(t,s)P(s)x\parallel^{p}
&=\int_{t-1}^{t} e^{p\gamma (t-s)}\parallel U(t,s)P(s)x\parallel^{p} d\tau\\
&\leq M^{p}e^{p\varepsilon t}\int_{t-1}^{t}e^{p(\gamma+\omega)(t- \tau)}e^{p\gamma(\tau-s)}\parallel U(\tau,s)P(s)x\parallel^{p} d\tau\\
&\leq M^{p}e^{p(\gamma+\omega)}e^{p\varepsilon t}\int_{s}^{\infty}e^{p\gamma(\tau-s)}\parallel U(\tau,s)P(s)x\parallel^{p} d\tau\\
&\leq K M^{p}e^{p(\gamma+\omega)}e^{p\varepsilon t}e^{p\delta s}\parallel x\parallel ^{p}
\end{align*}
and, respectively
\begin{align*}
e^{p\gamma (t-s)}\parallel V(s,t)Q(t)x&\parallel^{p}
=\int_{s}^{s+1}e^{p\gamma (t-s)}\parallel V(s,t)Q(t)x\parallel^{p} d\tau\\
&\leq M^{p}e^{p\varepsilon s}\int_{s}^{s+1}e^{p(\gamma+\omega)(\tau-s)}e^{p\gamma (t-\tau)}\parallel V(\tau,t)Q(t)x\parallel^{p} d\tau\\
&\leq M^{p}e^{p(\gamma+\omega)}e^{p\varepsilon s}\int_{0}^{t}e^{p\gamma (t-\tau)}\parallel V(\tau,t)Q(t)x\parallel^{p} d\tau\\
&\leq KM^{p}e^{p(\gamma+\omega)}e^{p\varepsilon s}e^{p\delta t}\parallel x\parallel ^{p}.
\end{align*}
Therefore,
\begin{equation}\label{De}
\parallel U(t,s)P(s)\parallel\leq K^{1/p}Me^{\gamma+\omega}e^{(\delta+\varepsilon) s}e^{-(\gamma-\varepsilon)(t-s)}
\end{equation}
and
\begin{equation}\label{Datko.2.1}
\parallel V(s,t)Q(t)\parallel\leq K^{1/p}Me^{\gamma+\omega}e^{(\delta+\varepsilon) t}e^{ -(\gamma+\varepsilon)(t-s)},
\end{equation}
for  $(t,s)\in\Delta_{+}$ with $t\geq s+1$.

If $t\in [s,s+1)$, it follows that
\begin{equation}\label{Da}
\parallel U(t,s)P(s)\parallel\leq Me^{\omega+\gamma-\varepsilon }e^{\varepsilon s}e^{-(\gamma-\varepsilon)(t-s)}
\end{equation}
and
\begin{equation}\label{Da.d}
\parallel V(s,t)Q(t) \parallel\leq Me^{\omega+\gamma}e^{\varepsilon t}e^{-(\gamma+\varepsilon)(t-s)}.
\end{equation}
By (\ref{De})--(\ref{Da.d}) we have  that there exist $N_1, N_2\geq 1$ such that
\begin{equation}\label{DD1}
\parallel U(t,s)P(s)\parallel\leq N_1e^{(\delta+\varepsilon) s}e^{-(\gamma-\varepsilon)(t-s)}
\end{equation}
and
\begin{equation}\label{DD2}
\parallel V(s,t)Q(t)\parallel\leq N_2e^{(\delta+\varepsilon) t}e^{ -(\gamma+\varepsilon)(t-s)},
\end{equation}
for all $(t,s)\in\Delta_{+}$.

Finally,  (\ref{DD1}) and (\ref{DD2}) involve that $U$ has a strong exponential dichotomy.
\end{proof}

\begin{remark}\label{rem.cond}\rm
If the evolution operator $U$ has a strong exponential
dichotomy with the dichotomy projection $P(\cdot)$ and constants $\alpha\geq 0$ and $\beta>0$ with $\alpha<\beta$, then  there exist constants $K\geq 1$, $\gamma>\varepsilon$ and $\delta\in[0,\gamma)$ such that $(\ref{eq.Datko})$ holds for each real parameter $p>0$. Indeed, in this case we can take $\varepsilon=\alpha$, $\delta=\alpha$ and $\gamma\in(\alpha,\beta)$.
Example \ref{ex.cond} shows that the assumptions $\gamma>\varepsilon$ and $\delta\in[0,\gamma)$ in the  theorem above are essential.
Furthermore, from the proof of the previous theorem (more precisely, from relations \eqref{DD1} and \eqref{DD2}) we observe that in order to have an exponential dichotomy (as in Definition \ref{n.e.d}), the assumption $\delta<\gamma$ can be omitted (we only need to assume $\gamma>\varepsilon$).
\end{remark}

In the uniform case we deduce the following result:

\begin{corollary}
Let  $U:\Delta_{+}\rightarrow\mathcal{B}(X)$ be an evolution operator and $P:\mathbb{R}_{+}\rightarrow \mathcal{B}(X)$ be a projection valued function compatible with $U$. The evolution operator $U$ has a uniform exponential dichotomy if and only if $U$ has $P$-uniform exponential growth and there is $p>0$ and there exist $K\geq 1$ and $\gamma>0$ such that
\begin{equation}\label{eq.Datko.uniform}
\int_{0}^{\infty}e^{p\gamma|\tau-t|}\parallel G(\tau,t)x\parallel ^{p}d\tau
\leq K|| x ||^p, \text{ for all  $t\geq 0$ and $x\in X$.}
\end{equation}
\end{corollary}
\begin{proof}
If $U$ has a uniform exponential dichotomy then, setting $\alpha=0$ in the proof of Proposition \ref{p1}, we have that $\varepsilon=\delta=0$ and thus inequality (\ref{eq.Datko.uniform}) holds for $\gamma\in(0,\beta)$ and $K=\max\left\{\frac{2N^p}{p(\beta-\gamma)},1\right\}$.
Conversely, letting $\varepsilon=\delta=0$ in relations (\ref{DD1}) and (\ref{DD2}) from Theorem \ref{p2}, we obtain that
$U$ has a uniform exponential dichotomy.
\end{proof}

\begin{remark}\label{D.uniform-}\rm
The previous result is also valid for $\gamma=0$ in inequality (\ref{eq.Datko.uniform}) (see \cite{Pr.Me.1985}).
\end{remark}


\subsection{A Lyapunov type theorem}
Given a real constant $\gamma>0$ and a projection valued function $P(\cdot)$, we denote by $\mathcal{H}_{\gamma}(P)$ the set of all strongly continuous operator valued functions
$ H:\mathbb{R}_{+}\rightarrow \mathcal{B}(X)$  with
\begin{equation}\label{eq.w}
\parallel H(t)x\parallel\leq e^{\gamma t}\parallel P(t)x\parallel+e^{-\gamma t}\parallel Q(t)x\parallel, \text{ for all } (t,x)\in\mathbb{R}_{+}\times X.
\end{equation}

Let $U:\Delta_{+}\rightarrow\mathcal{B}(X)$ be an evolution operator, $P:\mathbb{R}_{+}\rightarrow \mathcal{B}(X)$ be a projection valued function compatible with $U$ and let $H\in \mathcal{H}_{\gamma}(P)$. We say that a continuous function $L:\mathbb{R}_{+}\times X\longrightarrow \mathbb{R}$ is a \emph{Lyapunov function} corresponding to  $U$, $P(\cdot)$ and $H(\cdot)$ if
\begin{equation*}\label{eq.Ly}
L(t,U(t,s)x)+\int_{s}^{t}\parallel H(\tau)U(\tau,s)x\parallel^{2}d\tau\leq L(s,x), \text{ for all  $(t,s,x)\in\Delta_{+}\times X$.}
\end{equation*}

The following result shows that the existence of a strong exponential dichotomy implies the existence of a Lyapunov function.
\begin{theorem}\label{t.Ly.nec}
If  the evolution operator $U$ has a strong exponential
dichotomy with the dichotomy projection $P(\cdot)$ then  there exist $K\geq 1$, $\delta\geq 0$ and $\gamma>0$ such that for each $H\in \mathcal{H}_{\gamma}(P)$ there is a Lyapunov function $L$ which satisfies the following properties:
\begin{enumerate}

\item[$(L_1)$] $L(t,P(t)x)\geq 0$ and $L(t,Q(t)x)\leq 0$;

\item[$(L_2)$] $e^{-2\gamma t}L(t,P(t)x)-e^{2\gamma t}L(t,Q(t)x)\leq Ke^{2\delta t}\parallel x\parallel^{2}$;
\end{enumerate}
for all $(t,x)\in \mathbb{R}_{+}\times X$.
\end{theorem}
\begin{proof}
Setting $p=2$ in Proposition \ref{p1}, there exist  $K\geq 1$, $\delta\geq 0$ and $\gamma>0$ such that (\ref{eq.Datko}) holds.  For each  $H\in \mathcal{H}_{\gamma}(P)$,
taking
\begin{equation*}
L(t,x)=2\int_{0}^{\infty} sign({\tau-t}) \parallel H(\tau)G(\tau,t)x\parallel^{2}d\tau,
\end{equation*} we obtain

\begin{align*}
&L(t,U(t,s)x)+\int_{s}^{t}\parallel H(\tau)U(\tau,s)x\parallel^{2}d\tau=\\
&=2\int_{t}^{\infty} \parallel H(\tau)U(\tau,s)P(s)x\parallel^{2}d\tau-2\int_{0}^{s} \parallel H(\tau)V(\tau,s)Q(s)x\parallel^{2}d\tau-\\
&-2\int_{s}^{t} \parallel H(\tau)U(\tau,s)Q(s)x\parallel^{2}d\tau
+\int_{s}^{t}\parallel H(\tau)U(\tau,s)x \parallel^{2}d\tau\\
&\leq 2\int_{s}^{\infty} \parallel H(\tau)U(\tau,s)P(s)x\parallel^{2}d\tau-2\int_{0}^{s} \parallel H(\tau)V(\tau,s)Q(s)x\parallel^{2}d\tau\\
&=L(s,x), \text{ for  $(t,s,x)\in\Delta_{+}\times X$.}
\end{align*}
Therefore, $L$ is a Lyapunov function.
Moreover, for all  $(t,x)\in \mathbb{R}_{+}\times X$, we have
\begin{equation*}
L(t,P(t)x)=2\int_{t}^{\infty} \parallel H(\tau)U(\tau,t)P(t)x\parallel^{2}d\tau\geq 0
\end{equation*}
and, respectively
\begin{equation*}
L(t,Q(t)x)=-2\int_{0}^{t} \parallel H(\tau)V(\tau,t)Q(t)x\parallel^{2}d\tau\leq 0.
\end{equation*}
By (\ref{eq.w}) and Proposition \ref{p1}, we deduce that
\begin{align*}
&e^{-2\gamma t}L(t,P(t)x)-e^{2\gamma t}L(t,Q(t)x)=\\
&=2e^{-2\gamma t}\int_{t}^{\infty} \parallel H(\tau)U(\tau,t)P(t)x\parallel^{2}d\tau+2e^{2\gamma t}\int_{0}^{t} \parallel H(\tau)V(\tau,t)Q(t)x\parallel^{2}d\tau\\
&\leq 2 \int_{0}^{\infty}e^{2\gamma|\tau-t|}\parallel G(\tau,t)x\parallel ^{2}d\tau
\leq 2 Ke^{2\delta t}\parallel x\parallel^{2},
\end{align*}
for all  $(t,x)\in \mathbb{R}_{+}\times X$. This ends the proof.
\end{proof}

\begin{remark}\rm
The result in Theorem \ref{t.Ly.nec} remains valid even if the evolution operator $U$ has an exponential dichotomy. Furthermore,
by Remark \ref{rem.cond} and the proof of the theorem above, it follows that if the evolution operator $U$ has a strong exponential
dichotomy with $\alpha<\beta$, then  there exist a projection valued function $P(\cdot)$ compatible with $U$ and constants $K\geq 1$, $\gamma>\varepsilon$ and $\delta\in[0,\gamma)$ such that for each $H\in \mathcal{H}_{\gamma}(P)$ there is a Lyapunov function $L$ which satisfies $(L_1)$ and $(L_2)$.
\end{remark}

\begin{theorem}\label{t.Ly.suf}
Let $U:\Delta_{+}\rightarrow\mathcal{B}(X)$ be an evolution operator and $P(\cdot)$ be a projection valued function compatible with $U$.
If $U$ has $P$-exponential growth and there exist $K\geq 1$, $\gamma>\varepsilon$ and $\delta\in[0,\gamma)$ such that for each   $H\in \mathcal{H}_{\gamma}(P)$ there is a Lyapunov function $L$  which satisfies $(L_1)$ and $(L_2)$,
then $U$ has a strong exponential dichotomy.
\end{theorem}
\begin{proof} For  $\gamma>\varepsilon$, letting $$H(t)x=e^{\gamma t}P(t)x+e^{-\gamma t}Q(t)x, \text{ for } (t,x)\in\mathbb{R}_{+}\times X,$$
we have
\begin{align*}
&\int_{t}^{u}e^{2\gamma(\tau-t)}\parallel U(\tau,t)P(t)x\parallel ^{2}d\tau+\int_{0}^{t} e^{2\gamma(t-\tau)}\parallel V(\tau,t)Q(t)x\parallel^{2} d\tau\\
&= e^{-2\gamma t}\int_{t}^{u}\parallel H(\tau) U(\tau,t)P(t)x\parallel ^{2}d\tau+e^{2\gamma t}
\int_{0}^{t}\parallel H(\tau)V(\tau,t)Q(t)x\parallel ^{2}d\tau\\
&\leq e^{-2\gamma t}\left[ L(t,P(t)x)-L(u,U(u,t)P(t)x)\right]
+e^{2\gamma t}[L(0,V(0,t)Q(t)x)-L(t,Q(t)x)]\\
&\leq e^{-2\gamma t}L(t,P(t)x)-e^{2\gamma t}L(t,Q(t)x)
\leq K e^{2\delta t}\parallel x\parallel^{2},
\end{align*}
for all $u\geq t\geq 0$ and $x\in X$.
Letting $u\rightarrow\infty$, we obtain
\begin{align*}
\int_{0}^{\infty}e^{2\gamma|\tau-t|}\parallel G(\tau,t)x\parallel ^{2}d\tau
\leq K e^{2\delta t}\parallel x\parallel^{2}, \text{ for all  $(t,x)\in \mathbb{R}_{+}\times X.$}
\end{align*}
By Theorem \ref{p2}, we deduce that $U$ has a strong exponential dichotomy.
\end{proof}

\begin{example}\rm
Consider the evolution operator  $U$ and the projection valued function $P(\cdot)$ in Example \ref{ex.2} and take $K= 1$, $\delta=1$ and  $\gamma=2$.
For each strongly continuous operator valued function
$ H:\mathbb{R}_{+}\rightarrow \mathcal{B}(l^{\,\infty}(\mathbb{N},X))$  with
\begin{equation*}
\parallel H(t)x\parallel_{\infty}\leq e^{2 t}\sup\limits_{n\in\mathbb{N}} \parallel x_{2n}\parallel +e^{-2 t}\sup\limits_{n\in\mathbb{N}} \parallel x_{2n+1}\parallel, \text{ for  } (t,x)\in\mathbb{R}_{+}\times l^{\,\infty}(\mathbb{N},X)
\end{equation*}
(that is $H\in \mathcal{H}_{2}(P)$), we set
\begin{equation*}
L(t,x)=2\int_{t}^{\infty} \parallel H(\tau)(\frac{u(t)}{u(\tau)}x_0,0,\ldots)\parallel_{\infty}^{2}d\tau-2\int_{0}^{t} \parallel H(\tau)(0,\frac{u(\tau)}{u(t)}x_1,\ldots)\parallel_\infty^{2}d\tau,
\end{equation*}
for $x=(x_0,x_1,\ldots)\in l^{\,\infty}(\mathbb{N},X)$ and $t\geq 0$.
It follows from the proof of Theorem \ref{t.Ly.nec} that $L$ is a Lyapunov function. Moreover, conditions $(L_1)$ and $(L_2)$ hold. Thus, by Theorem \ref{t.Ly.suf}, we deduce that $U$ has a strong exponential dichotomy.
\end{example}

When $X$ is a  Hilbert space, we get the following result:
\begin{corollary}
Let $U:\Delta_{+}\rightarrow\mathcal{B}(X)$ be an evolution operator on the Hilbert space $X$ and $P:\mathbb{R}_{+}\rightarrow \mathcal{B}(X)$ be a projection valued function compatible with $U$.
If $U$ has  $P$-exponential growth and there exist $K\geq 1$, $\gamma>\varepsilon$ and $\delta\in[0,\gamma)$ such that for each   $H\in \mathcal{H}_{\gamma}(P)$ there is a strongly continuous operator valued function $W:\mathbb{R}_{+}\rightarrow \mathcal{B}(X)$ with   $W(t)^{*}=W(t),\,\forall\,t\geq 0$ and
\begin{enumerate}

\item[$(L_1')$] $\langle U(t,s)^{*}W(t)U(t,s)x+\int_{s}^{t}U(\tau,s)^{*}H(\tau)^{*}H(\tau)U(\tau,s)x\;d\tau,x\rangle\leq \langle W(s)x,x\rangle$;

\item[$(L_2')$] $\left<W(t)P(t)x,P(t)x\right>\geq 0$;

\item[$(L_3')$] $\left<W(t)Q(t)x,Q(t)x\right>\leq 0$;

\item[$(L_4')$] $e^{-2\gamma t}\left<W(t)P(t)x,P(t)x\right>-e^{2\gamma t}\left<W(t)Q(t)x,Q(t)x\right>\leq Ke^{2\delta t}\parallel x\parallel^{2}$;

\end{enumerate}
for all $(t,s,x)\in\Delta_+\times X$, then $U$ has a strong exponential dichotomy.
\end{corollary}
\begin{proof} It follows from Theorem \ref{t.Ly.suf}, letting
$$L(t,x)=\langle W(t)x,x\rangle,\text{ for all }(t,x)\in \mathbb{R}_{+}\times X. $$
\end{proof}

By Remark \ref{D.uniform-} and proceeding in a similar manner to the proofs of Theorem \ref{t.Ly.nec} and Theorem \ref{t.Ly.suf}, we obtain the following result for uniform exponential dichotomy:
\begin{corollary}
Let  $U:\Delta_{+}\rightarrow\mathcal{B}(X)$ be an evolution operator and $P:\mathbb{R}_{+}\rightarrow \mathcal{B}(X)$ be a projection valued function compatible with $U$. The evolution operator $U$ has a uniform exponential dichotomy if and only if $U$ has $P$-uniform exponential growth and there exist a function $L:\mathbb{R}_{+}\times X\longrightarrow \mathbb{R}$ and a constant $K\geq 1$ such that
\begin{enumerate}
\item  $L(t,U(t,s)x)+\int_{s}^{t}\parallel U(\tau,s)x\parallel^{2}d\tau\leq L(s,x)$;

\item  $L(t,P(t)x)\geq 0$;

\item  $L(t,Q(t)x)\leq 0$;

\item  $L(t,P(t)x)-L(t,Q(t)x)\leq K \parallel x\parallel^{2}$;
\end{enumerate}
for all $(t,s,x)\in\Delta_+\times X$.
\end{corollary}


\section{Open Problems}
An important property of exponential dichotomies is their roughness, that is, the preservation of exponential dichotomy
under sufficiently small perturbations of the evolution operator (we refer the reader to \cite{Ba.Va.2008-2,Ba.Va.2009-3,Po.2006} and the references therein).
We address the question of roughness for  strong exponential dichotomy.
Another open problem could be extending the results from \cite{Me.Sa.Sa.2003} to  strong exponential dichotomy.


\end{document}